\documentclass[10pt,letterpaper]{article}
\usepackage[top=0.85in,left=2.75in,footskip=0.75in]{geometry}

\usepackage{amsmath,amssymb}

\usepackage{changepage}

\usepackage[utf8x]{inputenc}

\usepackage{textcomp,marvosym}

\usepackage{cite}

\usepackage{nameref,hyperref}

\usepackage{microtype}
\DisableLigatures[f]{encoding = *, family = * }

\usepackage[table]{xcolor}

\usepackage{array}

\usepackage{algorithm, algorithmic}

\newcolumntype{+}{!{\vrule width 2pt}}

\newlength\savedwidth

\raggedright
\usepackage[aboveskip=1pt,labelfont=bf,labelsep=period,justification=raggedright,singlelinecheck=off]{caption}

\makeatletter
\renewcommand{\@biblabel}[1]{\quad#1.}
\makeatother

\usepackage{multirow}

\usepackage{lastpage,fancyhdr,graphicx}
\usepackage{epstopdf}
\fancyheadoffset[L]{2.25in}
\fancyfootoffset[L]{2.25in}
\begin{document}
\vspace*{0.2in}

\begin{flushleft}
{\Large
\textbf\newline{Application of Lie Group-based Neural Network Method to Nonlinear Dynamical Systems\footnote{Supported by  National Natural Science Foundation of China (11571008).}}  
}
\newline
\\
Ying Wen\textsuperscript{1\Yinyang},
Temuer Chaolu\textsuperscript{2*\Yinyang},
\\
\bigskip
\textbf{1} College of Information Engineering, Shanghai Maritime
University, 201306, Shanghai, China
\\
\textbf{2} College of Arts and Sciences, Shanghai Maritime
University, 201306, Shanghai, China
\\
\bigskip

%
%
\Yinyang These authors contributed equally to this work.





* tmchaolu@shmtu.edu.cn

\end{flushleft}
\section*{Abstract}
In this paper, a Lie group-based neural network method is proposed for solving initial value problems of non linear dynamics. Due to its single-layer structure (MLP), the approach is substantially cheaper than the multilayer perceptron method used in literature. The higher performance ability of the method is demonstrated by several examples.



\section{Introduction}
In recent years, a number of researchers have been devoted to the study of neural network (NN) approaches for solving differential equation (DE) problems. The study has been became a new developing trend of solving initial or boundary value (IVP or BVP) problem of an ordinary and partial differential equations (ODEs and PDEs). The main encourage of this study not only from the fact that NNs have general approximation ability \cite{Cybenko1989} but also due to their many advantages to study DEs, such as differ from many numerical methods to solve DEs, the methods do not rely on numeric difference scheme and can provide closed differentiable solutions to the problems of the DEs, etc.

Lee et al in \cite{Lee1990} solved the first order ODEs through using Hopfiled NN minimization methodology and compared the results to the traditional algorithm. Lagaris et al in \cite{Lagaris1998} proposed a NN trial solution method for the IVP of an ODE and PDEs; subsequently, irregular boundary problems were investigated by the method in \cite{Lagaris2000}; Maria et al in \cite{Piscopo2019} gave an extension on this method by training the network output directly as the equations solution, rather than relying on the trial solution, and successfully applied it to the calculation of the cosmological  phase transition tunneling profile. Baymani et al in \cite{Baymani2015} obtained the solutions of the Navier-Stokes equations in the form of analytical functions based on a NN method. Recently, researchers are moving further away from shallow NNs and toward deep learning. Raissi et al in \cite{Raissi2019} involves the development of data-driven solutions to PDEs and applies to quantum mechanics, reaction-diffusion systems, and so on. It is worthy to note that the network structure evolved as the research progressed, and new mathematical principles were incorporated. Chaouki et al in \cite{Aouiti2021} studied recurrent NNs with mixed time lag and time varying coefficients. Mall et al in \cite{Mall2014} employed Chebyshev networks to solve the singular IVP of Lane-Emden type equations, in which the computational difficulty is reduced through using Chebyshev polynomials with eliminate hidden layers in feedforward NNs. Marieme et al in \cite{Ngom2021} proposed a Fourier NN with the period function as the activation function, which simulated Fourier decomposition and was effectively applied to the Poisson equation and Heat equation. In a continuous-depth Bayesian NN family, Xu et al in \cite{Xu2021} conduct scalable approximate inference. In there, the hidden units follow a stochastic DE due to uncertainty about independent weights in each layer. To acquire continuous-time PDEs from sparse data, Iakovlev et al. in \cite{Iakovlev2020} used graphical NNs.

In most recent, as the extension and deepening of NN methods, people put forward various machine learning methods. However, the basis of these methods is still the forward NN algorithm (as the basic structural unit of these methods), which determines the ultimate efficiency of these algorithms.  At the same time, capturing the mathematical essence of the equation solution is the key to improve the efficiency of the algorithm. In real-world applications, these methods are also becoming more widely used. Eduardo et al in \cite{A2021} applied deep learning methods to solve PDEs for transport models. Zichao et al in \cite{PDE2017} proposed a new feed-forward deep network, called PDE-Net for estimating the dynamics of complex systems and their underlying implied PDE models. Most of the systems of DEs encountered in control applications are ODEs, and CK Ye et al in \cite{A2007} used artificial neural networks (ANNs) to solve the ODE in the state observers. Bing yu et al in \cite{Deep2019} applied deep learning techniques were applied to solve option prices using the backward stochastic differential equation (BSDE) method.

Although there are many researches, in considering convergence of the approximate solution, the most of literatures only tend to consider the equation itself with associated conditions, and do not sufficiently considered more implied information about the solution. This is one of the main reasons for some unsatisfactory performance of NN algorithm. For example, Lie algebraic structures, conservation laws, and possible expressions of solutions have not been applied in various machine learning methods. Raissia M et al in \cite{Raissi2019} proposed a physical information NN algorithm to solve DEs, in which some extra physical properties were tried to be used and got more accuracy solution to the problem. This enlightens more further study on this topics.

In this article, we try to use an algebra structure, admitted Lie groups, of underlying ODEs to get an alternative NN algorithm for solving the IVP of the nonlinear dynamic system.

A Lie group (symmetry) of a DEs is a one parameter transformations  mapping a solutio to another and providing formal exponential expression to the solution  of an IVP of the DEs. Lie group is vital to understanding the solutions characteristics of a nonlinear DEs since it provides a reasonable framework for analysing the solutions to the DEs. For further information on Lie groups of a DEs,  one refers to see \cite{Schwarz2007, Sharma2008, Tang2003, Devi2021}. Best in our knowledge, the essential combination of Lie group method, a powerful theoretical and computational tool, and NN method has not yet explored so far.

In current work, we propose a NN method with combining Lie group of a DEs to solve the IVP of the DEs. Unlike most cases, our proposed neural network algorithm with only one hidden layer can achieve high accuracy to approximate the solution. The simple structure of using a single hidden layer is beneficial not only in data preprocessing but subsequent stages of integration with any other more complex network. Our primary objective is to provide a general form of a NN solution for nonlinear dynamical systems by dividing the solution into two parts: the first term is a solution of an IVP of an ODEs, derived from the original IVP, can be solved by Lie group method without training; the second part is a NN with adjustable parameters, which will be learnt from an NN algorithm to make the final network solution approximates the solution of the DE. The NN is trained using an unsupervised approach throughout the procedure. The main advantages of the current method are: the method can provide more accurate solution expression, which will provide more accurate nonlinear information admitted by solution to make the subsequent network adjustment more effective; the effectiveness not only in the training interval, but also extrapolation outside the training interval. The efficiency and interpretability of the proposed algorithm are demonstrated by comparing the applications of the method with that of existing numerical method to some specific dynamics systems. Note that the idea used in this paper is a specific form of the Lie group method for the solution expression of the first-order ODEs rather that its reduction ability for a DEs (It can not reduce an IVP of a DEs). This is an innovative thinking and may provide some inspiration for further research on this issue.

The rest of present article is arranged as following. In Section II, we present the basic architecture of the Lie group-based NN and the learning algorithm. In Section III, we provide details of modeling examples and results as applications of our algorithm. In Section IV, we give some conclusion remarks.

\section{Preliminaries}
\subsection{System of initial value problems for first-order ordinary differential equations}
For a system of ODEs, the general form is given by
\begin{eqnarray}
\left\{\begin{array}{l}
y_{i}^{\prime}=f_{i}\left(t, y_{1}, y_{2}, \ldots, y_{n},\right) \\
y_{i}(0)=\alpha_{i},
\end{array} \quad 0 \leq t \leq a(i=0,1, \ldots, n) .\right. \label{ODEs}
\end{eqnarray}
where, $y_{i}(t)$ is the dependent variable, $t \in \mathbb{R}$ is the independent variable, and $f_{i}$ on the right side of the equation is a linear or nonlinear differentiable function about the independent and dependent variables.

We also know that first-order ODEs are just a special case of (\ref{ODEs}). Similarly, the IVPs for higher ODEs can be transformed into the forms (\ref{ODEs}).

\subsection{One-parameter continuous transformation group}
Suppose there are transformations $T_{\epsilon}$ in the $(x,y)$ plane:
\begin{eqnarray}
T_{\epsilon}: \hat{x}=\phi(x, y, \epsilon), \quad \hat{y}=\psi(x, y, \epsilon)  \label{transf}
\end{eqnarray}

\textbf{Definition 1}\ $G$ is a family of transformations $T_{\epsilon}$ which depends on the real parameters $\epsilon$ varying in an interval $\Delta$. If the following properties are satisfied, $G=\left\{T_{\epsilon}\right\}$ is called a one-parameter transformation group.

1) There are the following identity transformation when $\epsilon=0$.
\begin{eqnarray}
T_{0}: \hat{x}=\phi(x, y, 0), \quad \hat{y}=\psi(x, y, 0)  \label{e0}
\end{eqnarray}

2) When $\epsilon$ changes to $-\epsilon$, there is an inverse transformation.
\begin{eqnarray}
T_{\epsilon}: x=\phi(\hat{x}, \hat{y}, -\epsilon), \quad y=\psi(\hat{x},\hat{y}, -\epsilon)  \label{inverse}
\end{eqnarray}

3) Suppose
\begin{eqnarray}
T_{\delta}: x_{2}=\phi(\hat{x}, \hat{y}, \delta), \quad y_{2}=\psi(\hat{x}, \hat{y}, \delta) \label{product}
\end{eqnarray}
Then the product of these two transformations also belongs to (\ref{transf}). And the parameter becomes $\epsilon+\delta$, that is
\begin{eqnarray}
T_{\epsilon+\delta}: x_{2}=\phi(\hat{x}, \hat{y}, \delta)=\phi(x, y, \epsilon+\delta), \quad y_{2}=\psi(\hat{x}, \hat{y}, \delta)=\psi(x, y, \epsilon+\delta)  \label{new}
\end{eqnarray}

The $\phi(x, y, \epsilon)$, $\psi(x, y, \epsilon)$ in (\ref{transf}) is generally referred to as the global form of the continuous group. For a small parameter $\epsilon$, it can be expanded around $\epsilon=0$ as

\begin{eqnarray}
\hat{x}=x+\left.\epsilon \frac{\mathrm{d} \hat{x}}{\mathrm{~d} \epsilon}\right|_{\epsilon=0}+O\left(\epsilon^{2}\right), \quad \hat{y}=y+\left.\epsilon \frac{\mathrm{d} \hat{y}}{\mathrm{~d} \epsilon}\right|_{\epsilon=0}+O\left(\epsilon^{2}\right) \label{expand}
\end{eqnarray}

where, $O\left(\epsilon^{2}\right)$ represents the sum of all higher-order expansions including the second-order minima. Where the auxiliary functions $\xi(x, y)$ and $\eta(x, y)$ are introduced
\begin{eqnarray}
\xi=\left.\frac{\mathrm{d} \hat{x}}{\mathrm{~d} \epsilon}\right|_{\epsilon=0}, \quad \eta=\left.\frac{\mathrm{d} \hat{y}}{\mathrm{~d} \epsilon}\right|_{\epsilon=0} \label{xieta}
\end{eqnarray}
and relate it to the group $G$ given by (\ref{transf}), then we'll have
\begin{eqnarray}
\hat{x}=x+\epsilon \xi+O\left(\epsilon^{2}\right), \quad \hat{y}=y+\epsilon \eta+O\left(\epsilon^{2}\right)  \label{infinite}
\end{eqnarray}

This expansion of (\ref{expand}) around $\epsilon=0$ is called the infinitesimal of the transformation group (\ref{transf}).

\textbf{Theorem }\ The functions $\phi(x, y, \epsilon)$, $\psi(x, y, \epsilon)$ defining a group of transformations satisfy the system of differential equations
\begin{eqnarray}
\frac{d \hat{x}}{d \epsilon}=\xi(\hat{x}, \hat{y}), \quad \frac{d \hat{y}}{d \epsilon}=\eta(\hat{x}, \hat{y})  \label{gloform}
\end{eqnarray}
with initial value $\left.\hat{x}\right|_{\epsilon=0}=x$, $\left.\hat{y}\right|_{\epsilon=0}=y$. Conversely, when $\epsilon=0$ and the infinitesimal form of the continuous group is known, integration can be used to find the global form of the continuous group $G$.

\subsubsection{Infinitesimal operator of the group }
\textbf{Definition 2}\  An infinitesimal operator of the group $G$ is the linear differential operator
\begin{eqnarray}
X=\xi(x)^{i} \frac{\partial}{\partial x^{i}}  \label{x}
\end{eqnarray}
where $\xi(x)^{i}$ are determined in (\ref{xieta}). Functions $\xi(x)^{i}$ are coordinates of the operator $x$.

Let a function $f(x,y)$ be defined in the plane $(x,y)$, and under the transformation of (\ref{expand}), it becomes its overall similar $ f(\hat{x}, \hat{y})$. Expanding around $\epsilon=0$, as follows
\begin{eqnarray}
\begin{aligned}
f(\hat{x}, \hat{y}) &=f(x, y)+\epsilon X f+\frac{1}{2 !} \epsilon^{2} X^{2} f+\frac{1}{3 !} \epsilon^{3} X^{3} f+\cdots+\frac{1}{n !} \epsilon^{n} X^{n} f+\cdots \\
&=\sum_{n=0}^{\infty} \frac{1}{n !} \epsilon^{n} X^{n} f
\end{aligned}  \label{lieseries}
\end{eqnarray}

This series is called the Lie series. $X$ stands for Infinitesimal operator. From (\ref{x}), it is defined as
\begin{eqnarray}
X=\xi(x, y) \frac{\partial}{\partial x}+\eta(x, y) \frac{\partial}{\partial y}   \label{opex}
\end{eqnarray}
From (\ref{lieseries}) it can be seen that  Lie series can be expressed as
\begin{eqnarray}
f(\hat{x}, \hat{y})=\mathrm{e}^{\epsilon X} f(x, y)   \label{lies}
\end{eqnarray}

where $\mathrm{e}^{\epsilon X}=\sum_{n=0}^{\infty}(1 / n !) \epsilon^{n} X^{n}$ is an operator.
If we take $f(\hat{x},\hat{y})$ as $\hat{x}$ and $\hat{y}$ respectively in (\ref{lies}), we get
\begin{eqnarray}
\hat{x}=e^{\epsilon X} x, \quad \hat{y}=e^{\epsilon X} y  \label{lieses}
\end{eqnarray}

\subsection{Initial value problems for systems of first-order ordinary differential equations admitted by Lie group}
The IVPs of system of first-order ODEs is satisfied by the transformation $T_{\epsilon}: \hat{x}=\phi(x, y, \epsilon), \hat{y}=\psi(x, y, \epsilon)$ in the $G=\left\{T_{\epsilon}\right\}$ of (\ref{transf}). From (\ref{gloform}), there are the following
\begin{eqnarray}
\left\{\begin{aligned}
&\frac{d \hat{x}}{d \epsilon}=X \hat{x},\left.\hat{x}\right|_{\epsilon=0}=x\\
&\frac{d \hat{y}}{d \epsilon}=X \hat{y},\left.\hat{y}\right|_{\epsilon=0}=y
\end{aligned}\right. \label{lieode}
\end{eqnarray}

Under the transformation $\hat{x}=\phi(x, y, \epsilon), \hat{y}=\psi(x, y, \epsilon)$, the solution of (\ref{lieode}) are $\hat{x}=\mathrm{e}^{\epsilon X} x $ and $\hat{y}=\mathrm{e}^{\epsilon X} y $.

In \cite{Filippi1967} (The operator is denoted by $D$), the operator $X$ is divided into two main parts for convenience of calculation, $X=X_{1}+X_{2}$, and the solution of problem (\ref{lieode}) are rewritten as $\hat{x}=\mathrm{e}^{\epsilon (X_{1}+X_{2})} x $ and $\hat{y}=\mathrm{e}^{\epsilon (X_{1}+X_{2})} y $. We extend them to the following expressions
\begin{eqnarray}
\begin{aligned}
e^{\epsilon X} x=e^{\epsilon X_{1}} x+\sum_{\alpha=1}^{\infty} \sum_{k=\alpha}^{\infty} \frac{\epsilon^{k}}{k !} X_{1}^{k-\alpha} X_{2} X^{\alpha-1} x \\
e^{\epsilon X} y=e^{\epsilon X_{1}} y+\sum_{\alpha=1}^{\infty} \sum_{k=\alpha}^{\infty} \frac{\epsilon^{k}}{k !} X_{1}^{k-\alpha} X_{2} X^{\alpha-1} y  \label{X}
\end{aligned}
\end{eqnarray}
A proof of the $Gr\ddot{o}bner's$ formula for the Lie series and its convergence on (\ref{X}) is given in \cite{Filippi1967}.

To redescribe, the succinct equation has been devised. The Eq.(\ref{X}) are
\begin{eqnarray}
\left\{\begin{aligned}
&\hat{x}^{*}= e^{\epsilon X_{1}} x+\mathbb{N}_{1}(x,\epsilon) \\
&\hat{y}^{*}= e^{\epsilon X_{1}} y+\mathbb{N}_{2}(y,\epsilon)
\end{aligned} \right.\label{X1N}
\end{eqnarray}

From (\ref{lieode}), we know that $\bar{x}=e^{\epsilon X_{1}}x$ and $\bar{y}=e^{\epsilon X_{1}}y$ corresponding to the IVPs of a new system of ODEs $d\bar{x}/d\epsilon=X_{1} \bar{x}, d\bar{y}/d\epsilon=X_{1} \bar{y} $ with invariant initial value $\left.\bar{x}\right|_{\epsilon=0}=x$ and $\left.\bar{y}\right|_{\epsilon=0}=y$. The problem is simplified by this step. Solving a complex system of ODEs is extremely difficult, but choosing a simple component of the original equation makes the problem simple. In particular, it is also a partial solution of the original equation.

Some methods, such as Lie series, can be solved second portion $\mathbb{N}_{i}$ in (\ref{X1N}). When the equation is complex, it is difficult to solve. We propose that the second part be replaced with a black-box (a NN). It's worthy to note that $\mathbb{N}_{i}$ in (\ref{X1N}) are replaced by $\epsilon \mathbb{N}_{i}$ in order to satisfy the criterion that $\hat{x}^{*}(0)= x$ and $\hat{y}^{*}(0)= y$ when $\epsilon= 0$.

Based on the expression of the Lie group-based solution in (\ref{lieode}), Eq.(\ref{X1N}) are rewritten as
\begin{eqnarray}
\hat{x}^{*}= e^{\epsilon X_{1}} x+\epsilon \mathbb{N}_{1}(x,\epsilon), \quad \hat{y}^{*}= e^{\epsilon X_{1}} y+\epsilon \mathbb{N}_{2}(y,\epsilon) \label{X1XN}
\end{eqnarray}

The first part of the problem can be solved by select a simpler $X_{1}$, and the second part $\mathbb{N}_{i}$, should always be solved. The network solution $\hat{x}^{*}$ and $\hat{y}^{*}$ are being used to approximate the equation's solution $x$ and $y$. It is feasible to use the following expression
\begin{eqnarray}
\hat{x}\approx \bar{x}+\epsilon \mathbb{N}_{1}(\epsilon, p), \quad \hat{y} \approx \bar{y}+\epsilon \mathbb{N}_{2}(\epsilon, p) \label{X1NN}
\end{eqnarray}
The network output $\mathbb{N}_{i}$ is a value that incorporates just the input value $\epsilon$ and the network parameter $p$, with no equation boundary value.

The above is for the IVPs of a system of first-order ODEs, is also applicable to the IVPs of higher-order ODEs and first-order ODE. The mechanism is detailed below.

\section{Approximation of a system of ordinary differential equations by a neural network based on Lie groups}
In this section, we investigate the common types of nonlinear dynamical systems for which the Lie group-based NN algorithm solves the procedure as follows.

Consider the following coupled first-order ODEs,
\begin{eqnarray}
\left\{\begin{aligned}
&x^{\prime}=f_{1}\left(t, x, y\right), \\
&y^{\prime}=f_{2}\left(t, x, y\right).
\end{aligned}\right. \label{twoODEs}
\end{eqnarray}
With the initial value is $x(0)=\alpha_{0}$, $y(0)=\alpha_{1}$. convert it to the form (\ref{lieode}),
\begin{eqnarray}
\left\{\begin{aligned}
&\frac{d x}{d t}=Xx, \quad x(0)=\alpha_{0}\\
&\frac{d y}{d t}=Xy, \quad y(0)=\alpha_{1}
\end{aligned}\right. \label{twoodes}
\end{eqnarray}

The differential operator $X=f_{1}(t, x, y)\partial_{x}+f_{2}(t, x, y)\partial_{y}$. Written in the form of (\ref{X1XN}) as
\begin{eqnarray}
\left\{\begin{aligned}
&x=e^{tX}\alpha_{0},  \\
&y=e^{tX}\alpha_{1},
\end{aligned}\right. \label{twosol}
\end{eqnarray}

By the decomposability of the operator in (\ref{X}), Write the above equation in the form containing the network output
\begin{eqnarray}
\left\{\begin{aligned}
&\hat{x}=e^{tX_{1}}\alpha_{0}+t\mathbb{N}_{1}(t,p),  \\
&\hat{y}=e^{tX_{1}}\alpha_{1}+t\mathbb{N}_{2}(t,p),
\end{aligned}\right. \label{twonn}
\end{eqnarray}

Where $e^{tX_{1}}\alpha_{0}$ is the solution of equation $\bar{x}^{\prime}=X_{1}\bar{x}=\bar{f}_{1}(t,\bar{x},\bar{y})$ with initial value $\bar{x}(0)=x(0)=\alpha_{0}$, and $e^{tX_{1}}\alpha_{1}$ is the solution of equation $\bar{y}^{\prime}=X_{1}\bar{y}=\bar{f}_{2}(t,\bar{x},\bar{y})$ with initial value $\bar{y}(0)=y(0)=\alpha_{1}$. We clearly know that $e^{tX_{1}}\alpha_{0}$ is equivalent to $e^{t X_{1}} \alpha_{0}=\left.e^{t X_{1}} \bar{x}\right|_{\bar{x} \rightarrow \alpha_{0}}$, $e^{tX_{1}}\alpha_{1}$ is the same, and $t$ corresponds to the parameter $\epsilon$ of (\ref{lieode}). The number of equations determines $i$ in $\mathbb{N}_{i}$.

Any continuous function can be approximated by a feed-forward NN with a single hidden layer, according to the universal approximation theorem \cite{Cybenko1989}. This ANN ($\mathbb{N}$) can also be expressed as a matrix multiplication.
\begin{eqnarray}
N(t,p)=W_{2} \sigma\left(W_{1} t+b_{1}\right)+b_{2}   \label{NN}
\end{eqnarray}
$W_{1}\in R^{m \times n}$, which is the hidden layer weight from the input layer $t_{i}\in R^{n}(i=1,2,\ldots,n)$, $b_{1}\in R^{m}$, which is the hidden layer bias term, $W_{2}\in R^{m}$, which is the weight matrix of the weights from the hidden layer to the output layer, and $b_{2}\in R$, which is the output layer bias term. The number of neurons in the hidden layer is given by $m$, $\sigma=\tanh$ is a nonlinear activation function. The summation formula is expanded and written in vector form as $N(t, \vec{p})=\sum_{j=1}^{m}\left({w_{2}}_{j} \cdot \frac{e^{z_{j}}-e^{-z_{j}}}{e^{z_{j}}+e^{-z_{j}}}\right)+b_{2}$, $z_{j}=t \cdot {w_{1}}_{j}+{b_{1}}_{j}$, $j=1,2,\ldots,m$. The network structure is shown in Fig.\ref{fig1}. The number of networks grows in lockstep with the number of equations. $p$ represents the set of $W1, W2, b1, b2$. The NN takes a scalar input $t(t=\left\{t_{1},t_{2},\ldots,t_{n}\right\} in [0,a])$ and returns a scalar output $N(t)$, which fits the solution in (\ref{twoODEs}) of scalar functions $x(t)$ and $y(t)$ by $\hat{x}=\bar{x}+t N_{1}(t,p)$ and $\hat{y}=\bar{y}+t N_{2}(t,p)$. We could even use a more explicit expression afterwards.

\begin{figure}[!h]
\caption{{\bf Neural network models for systems of ODEs based on Lie groups}}
\label{fig1}
\end{figure}

\begin{eqnarray}
\hat{x}^{\prime}(t,p) \approx f_{1}(t,\hat{x}(t,p),\hat{y}(t,p)),\ \hat{y}^{\prime}(t,p) \approx f_{2}(t,\hat{x}(t,p),\hat{y}(t,p)) \label{fit}
\end{eqnarray}
This can be easily calculated analytically $x^{\prime}(t,p)$ and $y^{\prime}(t,p)$ without the need of techniques like differencing.
\begin{eqnarray}
\begin{aligned}
\hat{x}^{\prime}(t,p)=\frac{\partial\left(\bar{x}+t N_{1}(t,p)\right)}{\partial t}=\frac{\partial\bar{x}}{\partial t}+N_{1}(t,p)+t \frac{\partial N_{1}(t,p)}{\partial t}\\
\hat{y}^{\prime}(t,p)=\frac{\partial\left(\bar{y}+t N_{2}(t,p)\right)}{\partial t}=\frac{\partial\bar{y}}{\partial t}+N_{2}(t,p)+t \frac{\partial N_{2}(t,p)}{\partial t} \label{Diff}
\end{aligned}
\end{eqnarray}
$\frac{\partial\bar{x}}{\partial t}$ and $\frac{\partial\bar{y}}{\partial t}$ can be obtained by solving the differential equation by the properties of the Lie group.

The previously studied solution problem has been turned into an optimization problem in which the loss function $L(p)$ is minimized by optimizing the parameter $p$
\begin{eqnarray}
L(p)=\int_{0}^{a}\left[\hat{x}^{\prime}(t,p)-f_{1}(t,\hat{x}(t,p),\hat{y}(t,p))\right]^{2} d t\\
+\int_{0}^{a}\left[\hat{y}^{\prime}(t,p)-f_{2}(t,\hat{x}(t,p),\hat{y}(t,p))\right]^{2} d t \label{L}
\end{eqnarray}

To accomplish the calculation, we'll utilize the summation formula below.
\begin{eqnarray}
\begin{aligned}
L(p)&= L_{1}(p)+L_{2}(p) \\
&= \frac{1}{n}\left( {\sum_{i=1}^{n}\left[\hat{x}^{\prime}\left(t_{i},p\right)-f_{1}\left(t_{i},\hat{x}\left(t_{i},p\right),\hat{y}\left(t_{i},p\right) \right)\right]^{2} } \right. \\
 &\left.{+\sum_{i=1}^{n}\left[\hat{y}^{\prime}\left(t_{i},p\right)-f_{2}\left(t_{i},\hat{x}\left(t_{i},p\right),\hat{y}\left(t_{i},p\right) \right)\right]^{2}}\right) \\
\label{Loss}
 \end{aligned}
\end{eqnarray}
where $\left\{t_{i}\right\}$ is a set of training points that spans the $[0, a]$ domain. The number of $t$ values in the range $[0, a]$ is $n$. The error function $L_{1}$ relates to $\hat{x}$, while the error function $L_{2}$ relates to $\hat{y}$.

When the loss function $L(p)$ is small enough, the approach is capable of approximating the true solution $x(t)$ and $y(t)$ in the $[0, a]$. We not only use the above equation's mean square error to generate the loss function, but we also utilize the average root mean square error to assess the algorithm's superiority.
\begin{eqnarray}
L_{RMSE}=\frac{1}{2}\left(\mathcal{L}_{1}(p)+\mathcal{L}_{2}(p)\right)  \label{LRMSE}
\end{eqnarray}

where $\mathcal{L}_{1}(p)$ is the mean square error between the network solution $\hat{x}(t,p)$ and the exact solution $x(t)$(if the exact solution is unavailable, the numerical solution $x(t)$ is used), $\mathcal{L}_{1}(p)=\sqrt{\frac{1}{n}\left(\sum_{i=1}^{n}\left[\hat{x}\left(t_{i},p\right)-x\left(t_{i}\right)\right]^{2} \right)}$. Similarly, $\mathcal{L}_{2}(p)=\sqrt{\frac{1}{n}\left(\sum_{i=1}^{n}\left[\hat{y}\left(t_{i},p\right)-y\left(t_{i}\right)\right]^{2} \right)}$.

\textbf{Extension}\ This algorithm can also be used for IVPs for systems of higher-order ODEs, for the sake of brevity, will not be discussed here.

All the cases discussed in this paper are based on that scenario. The number of equations dictates the number of networks, and for two or more equations, the call between the solutions of the equations and the nonlinear characteristics of the equations causes the call between networks, and the final result is not as good as that of a single equation. To optimize network parameters and reduce the loss function $L(p)$, we used the more stable BFGS approach.

The proposed process approach is summarized in Algorithm \ref{alg:1}.

\begin{algorithm}
	\caption{The following are the steps in solving an ODEs using a feed forward NN model based on the Lie Group method.}
	\label{alg:1}
	\begin{algorithmic}[1]
		\STATE Consider the discrete points $t_{i}$ in the range $[0,a](i=1,2,...,n)$
		\STATE Configure the ANN architecture (number of layers, dimensionality of each layer and activation function, in this paper, we choose a hidden layer, the number of neurons in the hidden layer is $m$, and the activation function uses $\tanh$)
		\STATE Initialize the network parameters $p_{0}$, i.e., $W_{1}$, $b_{1}$, $W_{2}$, $b_{2}$
		\STATE Select the appropriate operator $X_{1}$ for (\ref{twoODEs})
		\STATE Calculate the objective function (trial solution) $\hat{x}=\bar{x}+tN_{1}(t, p)$ and $\hat{y}=\bar{y}+tN_{2}(t, p)$ by solving the new systems of equations corresponding to $X_{1}$
		\STATE Form the loss function $L(p)$ according to (\ref{Loss})
		\STATE To update the parameter $p$, adopt the BFGS optimization approach
        \STATE When the conditions are met, the optimal parameter $p^{*}$ is determined. The maximum number of iterations or the necessary accuracy of $L(p)$ can be selected
		\FORALL{$(W_{1}, b_{1}, W_{2}, b_{2}) \in p$}
		\STATE Update $W_{1}$ based on gradient $d L(p)/d W_{1}$ using  \textbf{$BFGS$}
		\STATE Update $b_{1}$ based on gradient $d L(p)/d b_{1}$ using \textbf{$BFGS$}
		\STATE Update $W_{2}$ based on gradient $d L(p)/d W_{2}$ using \textbf{$BFGS$}
		\STATE Update $b_{2}$ based on gradient $d L(p)/d b_{2}$ using \textbf{$BFGS$}
		\ENDFOR
		\STATE \textbf{return} $W_{1}$, $b_{1}$, $W_{2}$, $b_{2}$
        \STATE To determine the algorithm effectiveness, use the average RMSE (\ref{LRMSE})
	\end{algorithmic}
\end{algorithm}

\section{Numerical Experiments and Applications}
In this section, we perform experiments for some examples of nonlinear dynamical systems to demonstrate the behavior and properties of this new approach, comparing the accuracy and convergence of traditional numerical methods with the proposed method when the analytical solution is unknown, using the average mean square error to measure the accuracy of the numerical solution, and testing the algorithm's generalization capability. We select different types of practical problems and, where appropriate, point out the difficulties encountered. Although our NN just has one hidden layer, different models may use a different number of nodes for the length of the interval they examine. The following numerical experiments use the BFGS optimization algorithm (which is the most stable when compared to other optimization algorithms) to minimize the loss function in (\ref{Loss}), instead of directly calculating the derivatives of the NN, we use an automatic differentiation technique that has been proven to work well. Use the symbolic calculator mathematica to calculate the first part of the decomposition corresponding to the operator $X_{1}$. For the traditional numerical method using Scipy in Python.

\textbf{Example 1. } Ecosystems contain classic nonlinear dynamics equations for three-species food chains \cite{Global2003}, in which the bottom prey $x$ is preyed upon by the middle species $y$, which is then preyed upon by the top predator $z$. In the real world, the mouse-snake-owl is an example of a three-species food chain.
\begin{eqnarray}
\begin{aligned}
&\frac{d x}{d t}=a x-b x y
\end{aligned} \label{Lotka31}
\end{eqnarray}
\begin{eqnarray}
\begin{aligned}
&\frac{d y}{d t}=-c y+d x y-e y z
\end{aligned} \label{Lotka32}
\end{eqnarray}
\begin{eqnarray}
\begin{aligned}
&\frac{d z}{d t}=-f z+g y z
\end{aligned} \label{Lotka33}
\end{eqnarray}
where $x$ is the dependent variable representing the density of the bottom prey population, $y$ is the dependent variable representing the density of the middle predator(prey) population, $z$ is the dependent variable representing the density of the higher predator population, and $t$ is the independent variable representing time, which is $\epsilon$ in (\ref{X1NN}), the initial conditions are $x(0)=0.5, y(0)=1, z(0)=2$, and we consider the case where $t$ is in the interval $[0, 3]$. In the absence of a predator $(y = 0)$, Eq. (\ref{Lotka31}) shows that the prey $x$ would grow at a constant rate $a$, assuming that the prey $x$ have an unlimited food supply. Predation on the prey is considered to be proportional to the rate at which predators $y$ and prey $x$ are present at the same time, denoted by $bxy$. Predation is impossible if either $x$ or $y$ is zero. Similarly, in the absence of prey $(x = 0)$, Eq. (\ref{Lotka32}) demonstrates that the density of predators $y$ would decline at a constant rate $c$ due to natural death or emigration, $e$ represents the effect of predation on species $y$ by species $z$. This equation assumes that the predator population $y$ only hunts the same prey species as in Eq. (\ref{Lotka31}). In Eq. (\ref{Lotka33}), $f$ represents the natural death rate of species $z$ in the absence of prey, $g$ represents the reproduction rate of species $z$ in the presence of prey $y$. $a, b, c, d, e, f, g>0$.

According to the operator $X=\left(a x-b x y \right)\partial_{x} +\left(-c y+d x y-e y z\right)\partial_{y} + \left(-f z+g \right.$
$\left.y z \right)\partial_{z}$, from which the operator $X_{1}$ is chosen as $a x\partial_{x}-c y\partial{y}-f z\partial_{z}$, initial Parameters $a=b=c=d=e=f=g=1$, we can easily get the solution of this part $\bar{x}=\frac{1}{2} e^{t}, \bar{y}=e^{-t}, \bar{z}=2 e^{-t}$. Thus, the trial solution is $\hat{x}=\frac{1}{2} e^{t}+t N_{1}, \hat{y}=e^{-t}+t N_{2}, \hat{z}=2 e^{-t}+t N_{3}$.

We use a uniform distribution of $100$ training points in the interval $[0, 3]$ and $m=100$ to train the NN, similarly, $N_{2}$ and $N_{3}$ is a network model of this type. Since the model has no analytical solution, the fit to the numerical method is plotted as follows Fig.\ref{fig2}. The extrapolation capability of the test set at $[0, 3.5]$ is shown in Fig.\ref{fig3}. The optimization method is BFGS, as shown in Fig.\ref{fig4}, the log loss is the minimum compared to other methods at $1000$ iterations. The loss function of Eq.(\ref{Loss}), $L(p)=7.303\times10^{-5}$, $RMSE=0.00851$. We discovered during the experiment that when the training interval is more than $10$, the training effect does not reach the ideal state. It's also worth mentioning that the number of nodes necessary for training nonlinear differential equations in the hidden layer is more than before, which is another feature of nonlinear systems of equations.

\begin{figure}[!h]
\caption{{\bf Comparison between neural network solutions $\hat{x}$,$\hat{y}$,$\hat{z}$ and numerical solutions $x$,$y$,$z$ in Example 1}}
\label{fig2}
\end{figure}

\begin{figure}[!h]
\caption{{\bf Network network solution $\hat{x}$,$\hat{y}$,$\hat{z}$ and numerical solution $x$,$y$,$z$ in the test set}}
\label{fig3}
\end{figure}

\begin{figure}[!h]
\caption{{\bf Log loss function $L(p)$ in Example 1}}
\label{fig4}
\end{figure}

\textbf{Example 2. } We consider the van der Pol oscillator \cite{Period1987}, is a non-conservative oscillator with a linear spring force and a non-linear damping force. The equation is given
\begin{eqnarray}
\begin{aligned}
\frac{d^{2} x}{d t^{2}}-\mu\left(1-x^{2}\right) \frac{d x}{d t}+x=0
\end{aligned} \label{vander}
\end{eqnarray}

Applying the $Li\acute{e}nard$ transformation \cite{Transformation2000} $y=x-\frac{x^{3}}{3}-\frac{1}{\mu} \frac{d x}{d t}$,  the equation can be written as a system of ODEs,
\begin{eqnarray}
\left\{\begin{aligned}
&\frac{d x}{d t}=\mu\left(x-\frac{1}{3} x^{3}-y\right)\\
&\frac{d y}{d t}=\frac{x}{\mu}
\end{aligned}\right. \label{tranvander}
\end{eqnarray}

The operator form is $X=\mu\left(x-\frac{1}{3} x^{3}-y\right)\partial_{x}+\frac{x}{\mu}\partial_{y}$, We select to determine the analytical solution of $X_{1}=-\mu y \partial_{x}+\frac{x}{\mu}\partial_{y}$ from Eq.(\ref{tranvander}), $x(0)=1$, $y(0)=2$. The fit of the Runge-Kutta methods in the Scipy solver with the NN solution in the interval $[0,10]$, where $\mu=1$, is shown in Fig.\ref{fig5}. We use $40$ training data in the interval for training and $50$ nodes in the hidden layer. Fig.\ref{fig6} shows the extrapolation capability of the NN solution using trained parameters in the interval $[0, 11]$ compared to the numerical solution, again using BFGS optimization. Fig.\ref{fig7} displays the loss function iteration of the training, we can see that the error is about $10^{-1}$ after roughly $70$ iterations. Throughout the training process, the loss function $L(p)=2.07 \times 10^{-4}$, $RMSE=0.082$ is the measuring function.

\begin{figure}[!h]
\caption{{\bf Comparison between neural network solutions $\hat{x}$, $\hat{y}$ and numerical solutions $x$,$y$}}
\label{fig5}
\end{figure}

\begin{figure}[!h]
\caption{{\bf Network network solution $\hat{x}$, $\hat{y}$ and numerical solution $x$,$y$ in the Example 2}}
\label{fig6}
\end{figure}

\begin{figure}[!h]
\caption{{\bf Log loss function $L(p)$ in Example 2}}
\label{fig7}
\end{figure}

\textbf{Example 3. } A stimulation of chaos phenomena of Lorenz system. Consider the Lorenz equations \cite{Nonlinear1993}
\begin{eqnarray}
\left.\begin{array}{lll}
&\dot{y}_1=\sigma(y_2-y_1),\\
&\dot{y}_2=\rho y_1-y_2-y_1y_3,\\
&\dot{y}_3=-\beta y_3+y_1 y_2,
\end{array}
\right\}, (y_1,y_2,y_3)\in R^3,\sigma,\rho,\beta>0,\label{Lorenz}
\end{eqnarray}

As an alternative example of application of our method, we consider the case $\sigma=10,\rho=28, \beta=8/3$ and $y_1(0)=1, y_2(0)=5, y_3(0)=10$  discussed in \cite{Nonlinear1993}. The system has been already in the standard form. As previous examples, we take $X_1=10 y_2\partial_{y_1}-y_2 \partial_{y_2}-\beta y_3 \partial_{y_3}$. It yields the exact solutions to the associated Initial value problems of system (\ref{Lorenz})
\begin{eqnarray}
\left\{\begin{aligned}
&\bar{y}_1=e^{-10 t}\\
&\bar{y}_2=5 e^{- t}\\
&\bar{y}_3=10 e^{-28 t}\\
\end{aligned}\right. \label{LorenzD1}
\end{eqnarray}

Basing on the first parts of trial solutions $\hat{y}_i=\bar{y}_i+t N_i(t, p), i=1, 2, 3$. To derive the approximate solution $\hat{y}_{i}$, We trained the $N_{i}(t, p)$ network on the interval [0, 0.5] using uniformly spaced $40$ points with $30$ of hidden nodes. Fig.\ref{fig8} shows a comparison between the network solution and the numerical solution to (\ref{Lorenz}). Fig. \ref{fig9} shows that the ability to extrapolate in the interval $[0, 0.6]$ is not terrible. Fig. \ref{fig10} depicts the loss function convergence diagram. It can be seen that the optimization algorithm BFGS is the most stable and the easiest to reach the logarithmic minimum of the loss function.

\begin{figure}[!h]
\caption{{\bf Comparison between neural network solutions $\hat{y}_{i}$ and numerical solutions $y_{i}$}}
\label{fig8}
\end{figure}

\begin{figure}[!h]
\caption{{\bf Network network solution $\hat{y}_{i}$ and numerical solution $y_{i}$ in the Example 3}}
\label{fig9}
\end{figure}

\begin{figure}[!h]
\caption{{\bf Log loss function $L(p)$ in Example 3}}
\label{fig10}
\end{figure}

\textbf{Example 4.} We are thinking about another chaotic dynamical system consisting of three nonlinear ODEs, consider the following Rosseler system \cite{rossler1976}
\begin{eqnarray}
\begin{aligned}
&\dot{x}=-y-z \\
&\dot{y}=x+a y \\
&\dot{z}=b+z(x-c)
\end{aligned}\label{Rossler}
\end{eqnarray}

The operator $X_1=-z\partial_{x}+x\partial_{y}-c z\partial_{z}$ is determined from the operator $X=\left(-y-z\right)\partial_{x}+\left(x+a y\right)\partial_{y}+ \left(b+z(x-c)\right)\partial_{z}$, and this step is skipped owing to familiarity, the linear portion can be intercepted directly from Eq.(\ref{Rossler}) to produce the new systems of linear equations $\bar{x}^{\prime}=-\bar{z}$, $\bar{y}^{\prime}=\bar{x}$, and $\bar{z}^{\prime}=c \bar{z}$, with initial value $\bar{x}=1, \bar{y}=5, \bar{z}=10$. The trial solution is $\hat{x}=-0.754386+1.75439 e^{-5.7 t}+tN_{1}$, $\hat{y}=5.30779-0.307787 e^{-5.7 t}+tN_{2}$, and $\hat{z}=10 e^{-5.7 t}+tN_{3}$, as expected.

Our neural net solution achieves a close fit to the numerical method, as seen in Fig.\ref{fig11}, the hidden layer was trained with a grid of $40$ equidistant points in $t = [0, 1]$ and $50$ hidden units. $200$ points in the interval $[0, 1.4]$ are used to verify the trained network extrapolation capability, as shown Fig.\ref{fig12}. Fig.\ref{fig13} can be observed through a comparison of different optimization methods that the optimization method we employed is stable, loss function $L(p) = 3.266 \times 10^{-6}$ was achieved approximately $300$ times, which was used to estimate the algorithm $RMSE = 4.747 \times 10^{-5}$.

\begin{figure}[!h]
\caption{{\bf Comparison between neural network solutions $\hat{x}$,$\hat{y}$,$\hat{z}$ and numerical solutions $x$,$y$,$z$}}
\label{fig11}
\end{figure}

\begin{figure}[!h]
\caption{{\bf Network network solution $\hat{x}$,$\hat{y}$,$\hat{z}$ and numerical solution $x$,$y$,$z$ in the Example 4}}
\label{fig12}
\end{figure}

\begin{figure}[!h]
\caption{{\bf Log loss function $L(p)$ in Example 4}}
\label{fig13}
\end{figure}

\section{Discussion and conclusions}
Our network solves some useful problems and provides use cases for deep learning, and our network is structurally simpler and easier to study than other use cases. Our proposed lie group based NN algorithm,verified the accuracy of the method by solving the nonlinear model initial value problem. The initial network parameters are considered as random. From this study, it can be seen that the proposed model is easy to implement and the method is efficient and straightforward to solve for some nth ODEs or n sets of ODEs. It also provides a means to study nonlinear dynamical systems. We will utilize this method in PDEs in the future, or we will extend the network structure to deep networks, thinking about memory conditions and temporality through algorithmic upgrades.

\section*{Conflict of Interest:}
As far as we know, there are no conflicts of interest, financial or other conflicts between the designated author and the editors, reviewers and readers of this magazine.

\section*{Acknowledgements}
The authors thanks the supporting of National Natural Science Foundation of China with grand number 11571008.

%
%
%

\end{document}